\documentclass[a4paper,11pt]{amsart}

\usepackage[french]{babel}
\usepackage{amsmath}
\usepackage{amscd}
\usepackage{amssymb}
\usepackage[all]{xy}

\theoremstyle{plain}
\newtheorem{thm}{Th\'eor\`eme}[section]

\newtheorem{lem}[thm]{Lemme}
\newtheorem{prop}[thm]{Proposition}

\newtheorem*{thmintro}{Th\'eor\`eme}

\newtheorem*{propintro}{Proposition}

\theoremstyle{definition}

\theoremstyle{remark}

\newtheorem{rem}[thm]{Remarque}
\newtheorem{rems}[thm]{Remarques}

\newtheorem{question}[thm]{Question}

\swapnumbers
\theoremstyle{definition}
\newtheorem{setup}[thm]{Notations et hypoth\`eses}

\newcommand{\PP}{\textup{\textbf{P}}}
\newcommand{\CC}{\textup{\textbf{C}}}
\newcommand{\ZZ}{\textup{\textbf{Z}}}
\newcommand{\QQ}{\textup{\textbf{Q}}}

\renewcommand{\H}{\textup{H}}
\renewcommand{\P}{\textup{P}}
\renewcommand{\L}{\textup{L}}
\newcommand{\V}{\textup{V}}
\newcommand{\B}{\textup{B}}
\newcommand{\X}{\textup{X}}
\newcommand{\Y}{\textup{Y}}
\newcommand{\Z}{\textup{Z}}
\newcommand{\K}{\textup{K}}

\newcommand{\C}{\textup{C}}

\newcommand{\F}{\textup{F}}
\newcommand{\U}{\textup{U}}
\newcommand{\E}{\textup{E}}
\newcommand{\G}{\textup{G}}
\newcommand{\I}{\textup{I}}

\newcommand{\A}{\textup{A}}

\newcommand{\T}{\textup{T}}
\newcommand{\M}{\textup{M}}
\newcommand{\Q}{\textup{Q}}

\newcommand{\corps}{\textup{k}}

\renewcommand{\hom}{\textup{Hom}(\PP^{1},\textup{X})}
\newcommand{\homlibre}{\textup{Hom-l}(\PP^{1},\textup{X})}
\newcommand{\hombir}{\textup{Hom-b}(\PP^{1},\textup{X})}
\newcommand{\hombirn}{\textup{Hom-b}^{\textup{n}}(\PP^{1},\textup{X})}
\newcommand{\homx}{\textup{Hom}(\PP^{1},\textup{X},0\mapsto x)}
\newcommand{\homxlibre}{\textup{Hom-l}(\PP^{1},\textup{X},0\mapsto x)}
\newcommand{\homxbir}{\textup{Hom-b}(\PP^{1},\textup{X},0\mapsto x)}
\newcommand{\ux}{\textup{U}_{x}}
\newcommand{\vx}{\textup{V}_{x}}
\newcommand{\hx}{\textup{H}_{x}}
\newcommand{\secx}{\sigma_{x}}
\newcommand{\dx}{\textup{D}_{x}}

\newcommand{\ch}{\textup{ch}}
\newcommand{\td}{\textup{td}}

\newcommand{\groupe}{\textup{G}}
\newcommand{\stab}{\textup{G}_{0}}

\newcommand{\marque}{\stepcounter{thm}\noindent\textbf{\thethm. }}

\newcommand{\voir}{\textit{voir }}

\begin{document}

\title[Classes de Chern des vari\'et\'es unir\'egl\'ees]{Classes de Chern des vari\'et\'es unir\'egl\'ees}

\author{St\'ephane Druel}



\address{St\'ephane Druel, Institut Fourier, UMR 5582
du CNRS, Universit\'e Joseph Fourier, BP 74, 38402 Saint Martin d'H\`eres,
France.}

\email{druel@ujf-grenoble.fr}

%\urladdr{}

\maketitle

\thispagestyle{empty}

\section{Introduction}

Soient $\X$ une vari\'et\'e compacte de type k\"ahl\'erien et $\E$ est un fibr\'e vectoriel complexe de rang $r$ sur
$\X$. 

Le r\'esultat suivant montre qu'en g\'en\'eral peu de choses peuvent \^etre
dites de ses classes de Chern :
si $\X$ est une surface alg\'ebrique et si un entier $r>1$, $c_{1}\in\textup{NS}(\X)$ et
$c_{2}\in\textup{H}^{4}(\X,\ZZ)$ sont donn\'es alors il existe un fibr\'e vectoriel alg\'ebrique de rang $r$ sur $\X$
de classes de Chern $c_{1}$ et $c_{2}$.

Les exemples suivants montrent que sous certaines hypoth\`eses, il existe des contraintes sur les classes de Chern de
$\E$.
%Que peut-on dire et sous quelles hypoth\`eses de ses classes de Chern ?

Si $\E$ poss\`ede une m\'etrique d'Hermite-Einstein alors
l'in\'egalit\'e de Bogomolov-L\"ubke
$$\int_{\X}(2rc_{2}(\E)-(r-1)c_{1}(\E)^{2})\wedge\omega^{n-2}\ge 0$$
est satisfaite (\voir\cite{Lu82}). 
Si l'in\'egalit\'e pr\'ec\'edente est une \'egalit\'e alors $\E$ est projectivement plat. 

Si $\E$ est num\'eriquement effectif alors,
pour tout polyn\^ome de Schur $\P_{\lambda}\in\ZZ[c_{1}(\E),\ldots,c_{r}(\E)]$,
o\`u $\lambda$ est une partition d'un entier quelconque $k$ en parts $\le r$ (\voir\cite{DPS94} Theorem
2.5),
$$\int_{\X}\P(c(\E))\wedge\omega^{n-k}\ge 0.$$

Si $\X$ poss\`ede une m\'etrique de K\"ahler-Einstein et si $\omega$ est une forme de K\"ahler-Einstein alors sa
premi\`ere classe de Chern $c_{1}(\X)\in\H^{2}(\X,\CC)$ est n\'egative, nulle ou positive et l'in\'egalit\'e de
Miyaoka-Yau
$$\int_{\X}(2(n+1)c_{2}(\X)-nc_{1}(\X)^{2})\wedge\omega^{n-2}\ge 0$$
est satisfaite (\voir\cite{CO75}). Si l'in\'egalit\'e pr\'ec\'edente est une \'egalit\'e, alors le
rev\^etement universel de $\X$ est isomorphe \`a $\PP^{n}$, $\CC^{n}$ ou $\B^{n}$ (\voir\cite{Ti02} Theorem 2.13). 

\medskip
Nous d\'emontrons des r\'esultats d'effectivit\'e analogues au dernier de ceux que nous avons rappel\'es lorsque $\X$
est projective et unir\'egl\'ee.

\section{Enonc\'es des r\'esultats} \label{resultat}
\begin{setup}\label{setup}
Soit $\X$ une vari\'et\'e alg\'ebrique connexe, projective et lisse sur le corps $\CC$ des nombres complexes. Supposons
$\X$ unir\'egl\'ee, autrement dit, supposons que $\X$ soit domin\'ee par un produit $\PP^{1}\times\Y$ o\`u $\Y$ est une
vari\'et\'e alg\'ebrique complexe de dimension $\textup{dim}(\X)-1$.  

Soit $\H\subset\hombir$ une composante irr\'eductible telle que le morphisme d'\'evaluation
$\PP^{1}\times\H\longrightarrow\X$ soit dominant, autrement dit, \textit{une famille couvrante de courbes rationnelles}.
Soient $x\in\X$
et $\hx:=\homx\cap\H$. Soit $\stab\subset\textup{PGL(2)}$ le stabilisateur de $0\in\PP^{1}$. Le groupe $\stab$ 
op\`ere librement sur $\hx$ par la formule $g\cdot f=f\circ g^{-1}$
o\`u $g\in\stab$ et $[f]\in\hx$. Il op\`ere diagonalement sur $\PP^{1}\times\hx$. Soient
$\vx:=\hx//\stab$ et $\ux:=\PP^{1}\times\hx//\stab$ les quotients g\'eom\'etriques et $\pi_{x}$ et $\iota_{x}$
\begin{equation*}
\begin{CD}
\ux @){\iota_{x}})) \X \\
@V{\pi_{x}}VV\\
\vx
\end{CD}
\end{equation*}
les morphismes naturels. Supposons que pour $x\in\X$ g\'en\'eral, le
sch\'ema $\vx$ soit propre sur $\CC$, autrement dit, supposons que $\H$ soit \textit{une famille
couvrante minimale de courbes rationnelles.}
\end{setup}
\begin{rem}
Il suffit, par exemple, de consid\'erer une
famille couvrante de courbes 
rationnelles de degr\'e minimal relativement \`a une polarisation donn\'ee.
\end{rem}
Soit $\ell:=c_{1}(\X)\cdot f_{*}\PP^{1}\ge 2$ o\`u $[f]\in\H$. Les principaux r\'esultats de ce travail
sont les suivants.
\begin{thmintro}[\voir Th\'eor\`eme \ref{effectif}]
Sous les hypoth\`eses \ref{setup}, si $x\in\X$ est g\'en\'eral 
et si $\ux'$ est une composante irr\'eductible de $\ux$
alors le cycle 
$${\iota_{x}}_{*}(\ux')\cdot(c_{2}(\X)-\frac{\ell-1}{2\ell}c_{1}(\X)^{2})\in\A_{\ell-3}(\X)\otimes\QQ$$ 
est effectif.
\end{thmintro}
\begin{thmintro}[\voir Th\'eor\`eme \ref{principal}]
Sous les hypoth\`eses \ref{setup}, si $x\in\X$ est g\'en\'eral et si
$\ux'$ est une composante irr\'eductible de $\ux$
alors le cycle 
$${\iota_{x}}_{*}(\ux')\cdot(c_{2}(\X)-(\frac{\ell-1}{2\ell}+\frac{1}{\ell^{2}})c_{1}(\X)^{2})\in\A_{\ell-3}(\X)
\otimes\QQ$$ 
est effectif sauf s'il existe un morphisme fini 
$\widehat{\X}\longrightarrow\X$ et une application rationnelle
$\varphi:\widehat{\X}\dashrightarrow\Z$ dont les fibres g\'en\'erales sont des
espaces projectifs sur $\CC$ de dimension $\ell-1$, tels que les courbes rationnelles consid\'er\'ees soient les
images dans $\X$ des droites contenues dans les fibres de $\varphi$.
\end{thmintro}
\begin{propintro}[\voir Proposition \ref{egalite}]
Sous les hypoth\`eses \ref{setup}, le cycle
$${\iota_{x}}_{*}(\ux)\cdot(c_{2}(\X)-\frac{\ell-1}{2\ell}c_{1}(\X)^{2})\in\A_{\ell-3}(\X)\otimes\QQ$$ est nul
pour $x\in\X$ g\'en\'eral si et seulement s'il existe un morphisme fini 
$\widehat{\X}\longrightarrow\X$ et une application rationnelle
$\varphi:\widehat{\X}\dashrightarrow\Z$ dont les fibres g\'en\'erales sont des
espaces projectifs sur $\CC$ de dimension $\ell-1$, tels que les courbes rationnelles consid\'er\'ees soient les
images dans $\X$ des droites contenues dans les fibres de $\varphi$.
\end{propintro}
\begin{rems}
Si $\ell$ est un entier $\ge 2$ alors
$\displaystyle{c_{2}(\Q_{\ell-1})-(\frac{\ell-1}{2\ell}+\frac{1}{\ell^{2}})c_{1}(\Q_{\ell-1})^{2}=0}$ o\`u 
$\Q_{\ell-1}\subset\PP^{\ell}$ est une quadrique lisse de 
dimension $\ell-1$
et $\displaystyle{c_{2}(\PP^{\ell-1})-(\frac{\ell-1}{2\ell})c_{1}(\PP^{\ell-1})^{2}=0}$.

Si $\C$ est une courbe connexe, projective et lisse sur $\CC$ et si $\X_{\ell}:=\Q_{\ell-1}\times\C$
alors
$$c_{2}(\X_{\ell})-(\frac{\ell-1}{2\ell}+\frac{1}{\ell^{2}})c_{1}(\X_{\ell})^{2}
=\frac{\ell-1}{\ell^{2}}c_{1}(\Q_{\ell-1})\cdot c_{1}(\C).$$
Si le genre de $\C$ est $\ge 2$ et si $\ell\ge 3$ alors le cycle
$$c_{2}(\X_{\ell})-(\frac{\ell-1}{2\ell}+\frac{1}{\ell^{2}})c_{1}(\X_{\ell})^{2}$$
n'est pas effectif dans $\A_{\ell-3}(\X_{\ell})\otimes\QQ.$
\end{rems}
\begin{question}
Sous les hypoth\`eses \ref{setup}, si pour un point $x\in\X$ g\'en\'eral le
cycle
$${\iota_{x}}_{*}(\ux)\cdot(c_{2}(\X)-(\frac{\ell-1}{2\ell}+\frac{1}{\ell^{2}})c_{1}(\X)^{2})\in\A_{\ell-3}(\X)
\otimes\QQ$$ est nul, existe-t-il un morphisme fini $\widehat{\X}\longrightarrow\X$ et une application rationnelle
$\varphi:\widehat{\X}\dashrightarrow\Z$ dont les fibres g\'en\'erales sont des quadriques lisses sur $\CC$ de dimension
$\ell-1$,
tels que les courbes rationnelles consid\'er\'ees soient les
images dans $\X$ des droites contenues dans les fibres de $\varphi$ ?
\end{question}

Si $x\in\X$ est g\'en\'eral, $\vx$ est lisse sur $\CC$ de dimension
$\ell-2$ (\voir paragraphe \ref{tangent}) mais pas n\'ec\'essairement connexe. La description
globale du fibr\'e tangent de $\vx$ et le calcul de
$c_{1}(\vx)$ sont \'egalement donn\'es.
\begin{propintro}[\voir Proposition \ref{tangent quotient}]
Le fibr\'e tangent du sch\'ema $\vx$ est naturellement isomorphe au fibr\'e
$${\pi_{x}}_{*}((\iota_{x}^{*}\T_{X}/\T_{\pi_{x}})(-\secx)),$$ 
o\`u l'on a identifi\'e $\T_{\pi_{x}}\subset\T_{\ux}$ \`a son image dans $\iota_{x}^{*}\T_{\X}$ via la
diff\'erentielle de l'application $\iota_{x}$.
\end{propintro}
\begin{propintro}[\voir Proposition \ref{classe de chern}]
La premi\`ere classe de Chern de $\vx$ est donn\'ee par la formule
$$c_{1}(\vx)={\pi_{x}}_{*}\iota_{x}^{*}((\frac{\ell+1}{2\ell}-\frac{1}{\ell^{2}})c_{1}(\X)^{2}-c_{2}(\X))
\in\textup{Pic}(\vx)\otimes\QQ.$$
\end{propintro}
\section{Fibr\'e tangent du sch\'ema $\vx$}\label{tangent}
\marque Soit $\X$ une vari\'et\'e connexe, projective et lisse sur $\CC$. Soit $n$ la dimension
de $\X$ et supposons $n\ge 1$. Soient $\hom$ le sch\'ema localement quasi-projectif sur $\CC$ des morphismes de 
$\PP^{1}$ vers $\X$ et $\homx\subset\hom$ le sous-sch\'ema ferm\'e des morphismes
$f:\PP^{1}\longrightarrow\X$  tels que $f(0)=x$ o\`u $x\in\X$ est donn\'e 
(\voir\cite{Mo79} Proposition 1).
Soient 
$$\hombir\subset\hom$$ 
l'ouvert des morphismes birationnels sur leurs images et
$$\homxbir:=\homx\cap\hombir.$$
\noindent Un morphisme non constant $f\,:\,\PP^{1}\longrightarrow\X$ est dit libre si
$\textup{h}^{1}(\PP^{1},f^{*}\T_{\X}(-1))=0$, autrement dit, si
$\displaystyle{f^{*}\T_{\X}\simeq\mathcal{O}_{\PP^{1}}(a_{1})\oplus\cdots\oplus\mathcal{O}_{\PP^{1}}(a_{n})}$
avec $a_{1}\ge\cdots\ge a_{n}\ge 0$ (\voir\cite{Ko96} Definition II 3.1). 
Soient enfin 
$$\homlibre\subset\hom$$ 
l'ouvert des morphismes libres (\voir\cite{Ko96} Corollary II 3.5.4) et
$$\homxlibre:=\homx\cap\homlibre.$$
\noindent Le sch\'ema $\homxlibre$ est lisse sur $\CC$ de dimension $\textup{h}^{0}(\PP^{1},f^{*}\T_{\X}\otimes\I_{0})$
en
$[f]\in\homxlibre$ (\voir\cite{Mo79}  Proposition 2). 

\medskip

\marque Supposons $\X$ unir\'egl\'ee. Soient $\H\subset\hombir$ une famille couvrante de courbes rationnelles et
$\hx:=\homx\cap\H$ o\`u $x\in\X$ est donn\'e. Soient $ev_{x}$ et $q_{x}$
\begin{equation*}
\begin{CD}
\PP^{1}\times\hx @){ev_{x}})) \X \\
@V{q_{x}}VV\\
\hx
\end{CD}
\end{equation*}
les morphismes naturels.
\begin{lem}[\voir\cite{De01} Proposition 4.14]\label{libre}
Si $x\in\X$ est g\'en\'eral alors les morphismes param\'etr\'es par $\hx$
sont libres.
\end{lem}
Le sch\'ema $\hx$ est en particulier lisse sur $\CC$ de dimension $\textup{h}^{0}(\PP^{1},f^{*}\T_{\X}\otimes\I_{0})$
en $[f]\in\hx$.
\begin{prop}\label{tangent hom}
Si $x\in\X$ est g\'en\'eral alors le fibr\'e tangent du sch\'ema $\hx$ est naturellement isomorphe au fibr\'e
$${q_{x}}_{*}(ev_{x}^{*}\T_{X}(-\{0\}\times\hx).$$
\end{prop}
\begin{proof}
La restriction de la diff\'erentielle de $ev_{x}$ \`a $q_{x}^{*}\T_{\hx}$ est 
identiquement nulle le long de $\{0\}\times\hx\subset\PP^{1}\times\hx$ et d\'efinit donc
un morphisme
$$q_{x}^{*}\T_{\hx}\longrightarrow ev_{x}^{*}\T_{X}(-\{0\}\times\hx)$$
et, en prenant les images directes par $q_{x}$, un morphisme 
$$\T_{\hx}\longrightarrow {q_{x}}_{*}(ev_{x}^{*}\T_{X}(-\{0\}\times\hx))$$
dont nous allons montrer qu'il est l'isomorphisme cherch\'e.

L'espace tangent de Zariski \`a $\hx$ en $[f]\in\hx$
s'identifie naturellement \`a $\textup{H}^{0}(\PP^{1},f^{*}\T_{\X}\otimes\I_{0})$ (\voir\cite{Mo79} Proposition 2)
et la diff\'erentielle de l'application $ev_{x}$ en 
$(t,f)\in\PP^{1}\times\hx$ (\voir\cite{Ko96} Proposition II 3.4)
est l'application 
$$\T_{\PP^{1}}\otimes\corps(t)\oplus
\textup{H}^{0}(\PP^{1},f^{*}\T_{\X}\otimes\I_{0})\longrightarrow\T_{\X}\otimes\corps(f(t))\simeq(f^{*}\T_{\X})
\otimes\corps( t )$$
donn\'ee par
$$(v,s)\mapsto \textup{d}f_{t}(v)+s(t),$$
\noindent o\`u $v\in\T_{\PP^{1}}\otimes\corps(t)$ et
$s\in\textup{H}^{0}(\PP^{1},f^{*}\T_{\X}\otimes\I_{0})$. La compos\'ee de l'application 
$$\textup{H}^{0}(\PP^{1},f^{*}\T_{\X}\otimes\I_{0})\simeq\T_{\hx}\otimes\corps(f)\longrightarrow
({q_{x}}_{*}(ev_{x}^{*}\T_{X}(-\{0\}\times\hx))\otimes\corps(f)$$ 
et de l'application naturelle d'\'evaluation des sections  
$$({q_{x}}_{*}(ev_{x}^{*}\T_{X}(-\{0\}\times\hx))\otimes\corps(f)\longrightarrow 
\textup{H}^{0}(\PP^{1},f^{*}\T_{\X}\otimes\I_{0})$$ 
est donc l'identit\'e de $\textup{H}^{0}(\PP^{1},f^{*}\T_{\X}\otimes\I_{0})$.

Il reste \`a remarquer que la formation de l'image directe
${q_{x}}_{*}(ev_{x}^{*}\T_{X}(-\{0\}\times\hx)$
commute aux changements de bases, autrement dit, que l'application naturelle 
$$({q_{x}}_{*}(ev_{x}^{*}\T_{X}(-\{0\}\times\hx))\otimes\corps(f)\longrightarrow
\textup{H}^{0}(\PP^{1},f^{*}\T_{\X}\otimes\I_{0})$$
est un isomorphisme (\voir\cite{Ha77} III Theorem 12.11).
\end{proof}
\begin{rems}
L'\'enonc\'e pr\'ec\'edent reste vrai si $x\in\X$ est
quelconque et si $\hx$ est une composante connexe de $\homxlibre$.

Soit $\H$ une composante connexe de $\homlibre$ et soient $ev$ et $q$
\begin{equation*}
\begin{CD}
\PP^{1}\times\H @){ev})) \X \\
@V{q}VV\\
\H
\end{CD}
\end{equation*}
les morphismes naturels. La preuve de la proposition suivante est analogue \`a celle de la proposition
\ref{tangent hom}.
\begin{prop}
Le fibr\'e tangent du sch\'ema $\H$ est naturellement isomorphe au fibr\'e $q_{*}ev^{*}\T_{X}.$
\end{prop}
\end{rems}
\marque Soient $\groupe:=\textup{PGL(2)}$ et
$\stab\subset\groupe$ le stabilisateur de $0\in\PP^{1}$. 
Le groupe $\stab$ op\`ere librement sur $\hx$
par la formule $g\cdot f=f\circ g^{-1}$
o\`u $g\in\stab$ et $[f]\in\hx$ (\voir\cite{Mo79} Lemma 9). 
Il op\`ere diagonalement sur $\PP^{1}\times\hx$. Les quotients g\'eom\'etriques
$$\vx:=\hx//\stab\quad\text{et}\quad\ux:=\PP^{1}\times\hx//\stab$$ 
existent et 
les vari\'et\'es $\hx$ et $\PP^{1}\times\hx$ 
sont des fibr\'es principaux de groupe $\stab$ au dessus des quotients $\vx$ et $\ux$ respectivement
(\voir\cite{Mo79} Lemma 9). Les vari\'et\'es $\vx$ et $\ux$ sont en particulier lisses sur $\CC$ si $x\in\X$ est
g\'en\'eral. Soient $\pi_{x}$ et $\iota_{x}$
\begin{equation*}
\begin{CD}
\ux @){\iota_{x}})) \X \\
@V{\pi_{x}}VV\\
\vx
\end{CD}
\end{equation*}
les morphismes naturels : $\pi_{x}$ est une fibration
en droites projectives sur $\CC$. Soit $\secx\subset\ux$ la section de $\pi_{x}$, quotient g\'eom\'etrique de
$\{0\}\times\hx$ par le groupe $\stab$. La section $\secx$ est contract\'ee par $\iota_{x}$ sur $x\in\X$.
\begin{prop}\label{tangent quotient}
Si $x\in\X$ est g\'en\'eral alors le fibr\'e tangent du sch\'ema $\vx$ est naturellement isomorphe au fibr\'e
$${\pi_{x}}_{*}((\iota_{x}^{*}\T_{X}/\T_{\pi_{x}})(-\secx)),$$ 
o\`u l'on a identifi\'e $\T_{\pi_{x}}\subset\T_{\ux}$ \`a son image dans $\iota_{x}^{*}\T_{\X}$ via la
diff\'erentielle de l'application $\iota_{x}$.
\end{prop}
\begin{proof}
Notons $r_{x}$ et $s_{x}$

\centerline{
\xymatrix{
\PP^{1}\times\hx \ar[r]^-{r_{x}} \ar[d]_-{q_{x}} \ar@/^1.5pc/[rr]^-{ev_{x}}
& \ux \ar[r]^-{\iota_{x}} \ar[d]_-{\pi_{x}} & \X \\
\hx \ar[r]^-{s_{x}} & \vx
}
}
\noindent les morphismes de passage au quotient. Le fibr\'e vectoriel $s_{x}^{*}\T_{\vx}$ s'identifie,
via
la diff\'erentielle de $s_{x}$, au conoyau de l'application
$$\T_{s_{x}}\longrightarrow\T_{\hx}.$$
\indent La diff\'erentielle de l'action de $\stab$ sur $\hx$ en
$(e,f)\in\stab\times\hx$ est l'application 
$$\textup{H}^{0}(\PP^{1},\T_{\PP^{1}}\otimes\I_{0})\oplus\textup{H}^{0}(\PP^{1},f^{*}\T_{\X}\otimes\I_{0})
\longrightarrow\textup{H}^{0}(\PP^{1},f^{*}\T_{\X}\otimes\I_{0})$$
donn\'ee par
$$(x,y)\mapsto y-f_{*}x,$$
o\`u $x\in\textup{H}^{0}(\PP^{1},\T_{\PP^{1}}\otimes\I_{0})$ et
$y\in\textup{H}^{0}(\PP^{1},f^{*}\T_{\X}\otimes\I_{0})$.
La formation des images directes ${q_{x}}_{*}\T_{q_{x}}(-\{0\}\times\hx)$ et
${q_{x}}_{*}(ev_{x}^{*}\T_{X}(-\{0\}\times\hx)$ 
commute aux changements de bases, autrement dit, les applications naturelles
$$({q_{x}}_{*}\T_{q_{x}}(-\{0\}\times\hx))\otimes\corps(f)\longrightarrow
\textup{H}^{0}(\PP^{1},\T_{\PP^{1}}\otimes\I_{0})$$
et 
$$({q_{x}}_{*}(ev_{x}^{*}\T_{X}(-\{0\}\times\hx))\otimes\corps(f)\longrightarrow
\textup{H}^{0}(\PP^{1},f^{*}\T_{\X}\otimes\I_{0})$$
sont des isomorphismes (\voir\cite{Ha77} III Theorem 12.11). Le sous-fibr\'e vectoriel
$$\T_{s_{x}}\subset\T_{\hx}$$ 
s'identifie donc au sous-fibr\'e vectoriel
$${q_{x}}_{*}\T_{q_{x}}(-\{0\}\times\hx)\subset{q_{x}}_{*}(ev_{x}^{*}\T_{X}(-\{0\}\times\hx),$$
via l'identification $\T_{\hx}\simeq{q_{x}}_{*}(ev_{x}^{*}\T_{X}(-\{0\}\times\hx)$ donn\'ee par la proposition
\ref{tangent hom}.
Ces identifications sont $\stab$-\'equivariantes.

Le diagramme 
\begin{equation*}
\begin{CD}
\PP^{1}\times\hx @){r_{x}})) \ux \\
@V{q_{x}}VV  @V{\pi_{x}}VV \\
\hx @){s_{x}})) \vx
\end{CD}
\end{equation*}
est un carr\'e cart\'esien dans la cat\'egorie des $\stab$-sch\'emas puisque l'action de $\stab$ sur $\hx$ est libre. 
La restriction de la diff\'erentielle de $r_{x}$ \`a $\T_{q_{x}}$ induit un isomorphisme \'equivariant de $\T_{q_{x}}$
sur $r_{x}^{*}\T_{\pi_{x}}$.
Le morphisme $s_{x}$ \'etant plat, l'application
$${q_{x}}_{*}r_{x}^{*}(\T_{\pi_{x}}(-\secx))\simeq{q_{x}}_{*}\T_{q_{x}}(-\{0\}\times\hx)
\longrightarrow{q_{x}}_{*}(ev_{x}^{*}\T_{X}(-\{0\}\times\hx)\simeq{q_{x}}_{*}r_{x}^{*}(\iota_{x}^{*}\T_{\X}(-\secx))$$
s'identifie donc \`a l'application
$$s_{x}^{*}{\pi_{x}}_{*}\T_{\pi_{x}}(-\secx)\longrightarrow s_{x}^{*}{\pi_{x}}_{*}ev_{x}^{*}\T_{X}(-\secx)$$
de fa\c{c}on \'equivariante.
L'isomorphisme \'equivariant
$$s_{x}^{*}\T_{\vx}\simeq s_{x}^{*}{\pi_{x}}_{*}((\iota_{x}^{*}\T_{X}/\T_{\pi_{x}})(-\secx))$$ 
descend en l'isomorphisme cherch\'e
$$\T_{\vx}\simeq{\pi_{x}}_{*}((\iota_{x}^{*}\T_{X}/\T_{\pi_{x}})(-\secx)).$$
\end{proof}
\begin{rems}
L'\'enonc\'e pr\'ec\'edent reste vrai si $x\in\X$ est quelconque et si $\hx$ est une
composante connexe de $\homxlibre\cap\hombir$.

Soit $\H$ une composante connexe de $\homlibre\cap\hombir$.
Le groupe $\groupe$ op\`ere librement sur $\H$
par la formule
$g\cdot f=f\circ g^{-1}$
o\`u $g\in\groupe$ et $[f]\in\H$ (\voir\cite{Mo79} Lemma 9). 
Il op\`ere diagonalement sur $\PP^{1}\times\H$. Les quotients g\'eom\'etriques
$\V:=\H//\groupe$ et 
$\U:=\PP^{1}\times\H//\groupe$
existent et les sch\'emas $\H$ et $\PP^{1}\times\H$ 
sont des fibr\'es principaux de groupe $\groupe$ au dessus des quotients $\V$ et $\U$ respectivement
(\voir\cite{Mo79} Lemma 9). Les sch\'emas $\V$ et $\U$ sont en particulier lisses sur $\CC$. Soient
$\pi$ et $\iota$
\begin{equation*}
\begin{CD}
\U @){\iota})) \X \\
@V{\pi}VV\\
\V
\end{CD}
\end{equation*}
les morphismes naturels. Le morphisme $\pi$ est une fibration
en droites projectives sur $\CC$. La preuve de la proposition suivante est analogue \`a celle de la proposition
\ref{tangent quotient}.
\begin{prop}
Le fibr\'e tangent du sch\'ema $\V$ est naturellement isomorphe au fibr\'e
$\pi_{*}(\iota^{*}\T_{X}/\T_{\pi}),$
o\`u l'on a identifi\'e $\T_{\pi}\subset\T_{\U}$ \`a son image dans $\iota^{*}\T_{\X}$ via la
diff\'erentielle de l'application $\iota$.
\end{prop}
\end{rems}
\section{Un calcul de classe de Chern}\label{calcul}
\marque Nous reprenons les notations de la section pr\'ec\'edente. Soit $x\in\X$ un point g\'en\'eral. 
Soit $\E_{x}$ le fibr\'e vectoriel de rang 2 sur $\vx$
donn\'e par $\E_{x}:={\pi_{x}}_{*}\mathcal{O}_{\ux}(\secx)$. La section $\secx$ de $\pi_{x}$ correspond \`a la donn\'ee
d'un quotient inversible $\M_{x}$ de $\E_{x}$. Le fibr\'e $\E_{x}$ est une extension 
$$0\longrightarrow\mathcal{O}_{\vx}\longrightarrow\E_{x}\longrightarrow\M_{x}\longrightarrow 0$$
du fibr\'e inversible $\M_{x}$ par le fibr\'e trivial $\mathcal{O}_{\vx}$ et $c_{2}(\E_{x})=0$. Soit
$\L_{x}:=\M_{x}^{\otimes -1}$.
Le sch\'ema $\ux$ s'identifie au fibr\'e projectif $\PP_{\vx}(\E_{x})$ au dessus de $\vx$ et
l'id\'eal de $\secx$ dans $\ux$ s'identifie au fibr\'e tautologique $\mathcal{O}_{\PP_{\vx}(\E_{x})}(-1)$.

Soit $[f]\in\hx$ et soit $\ell:=c_{1}(\X)\cdot f_{*}(\PP^{1})\ge 2$. Le fibr\'e en droites
$\iota_{x}^{*}\omega_{\X}$ est de degr\'e $-\ell$ sur les fibres de $\pi_{x}$ et
trivial
sur $\secx$. Il est donc isomorphe au produit de $\ell$ copies du fibr\'e $\mathcal{O}_{\PP_{\vx}(\E_{x})}(-1)\otimes
\pi_{x}^{*}\M_{x}$.
\begin{prop}\label{classe de chern}
Si $x\in\X$ est g\'en\'eral alors la premi\`ere classe de Chern de $\vx$ est donn\'ee par la formule
$$c_{1}(\vx)={\pi_{x}}_{*}\iota_{x}^{*}((\frac{\ell+1}{2\ell}-\frac{1}{\ell^{2}})c_{1}(\X)^{2}-c_{2}(\X))
\in\textup{Pic}(\vx)\otimes\QQ.$$
\end{prop}
\begin{proof}
Soient $\textup{A}(\vx)$ l'anneau de Chow de $\vx$ et $\textup{K}(\vx)$ son anneau de Grothendieck.
Soient $\ch\,:\,\textup{K}(\vx)\longrightarrow\textup{A}(\vx)$ le caract\`ere de Chern exponentiel et
$\td(\T_{\pi_{x}})\in\textup{A}(\ux)$ 
la classe de Todd du fibr\'e tangent $\T_{\pi_{x}}$ \`a $\pi_{x}$. La formule utile est 
la formule de Grothendieck-Riemann-Roch (\voir\cite{Fu98} Theorem 15.2) :
$$\ch({\pi_{x}}_{!}((\iota_{x}^{*}\T_{X}/\T_{\pi_{x}})(-\secx)))
={\pi_{x}}_{*}(\ch((\iota_{x}^{*}\T_{X}/\T_{\pi_{x}})(-\secx))
\cdot \td (\T_ {\pi_{x}}))\in\textup{A}(\vx)\otimes\QQ.$$
Les morphismes param\'etr\'es par $\hx$ sont libres d'apr\`es le lemme \ref{libre}.
Le faisceau 
$$\textup{R}^{1}{\pi_{x}}_{*}((\iota_{x}^{*}\T_{X}/\T_{\pi_{x}})(-\secx))$$
est donc nul,
$${\pi_{x}}_{!}((\iota_{x}^{*}\T_{X}/\T_{\pi_{x}})(-\secx))={\pi_{x}}_{*}((\iota_{x}^{*}\T_{X}/\T_{\pi_{x}})(-\secx))$$
et, d'apr\`es la proposition \ref{tangent quotient}, 
$$\ch({\pi_{x}}_{!}((\iota_{x}^{*}\T_{X}/\T_{\pi_{x}})(-\secx)))=\ch(\T_{\vx}).$$
La premi\`ere classe de Chern de $\vx$ est donc donn\'ee par la formule
$$c_{1}(\vx)=\ch(\T_{\vx})_{1}={\pi_{x}}_{*}(\ch((\iota_{x}^{*}\T_{X}/\T_{\pi_{x}})(-\secx))\cdot\td(\T_{\pi_{x}}))_{2}.
$$
Il reste \`a calculer la composante de codimension 2 du cycle
$$\ch((\iota_{x}^{*}\T_{X}/\T_{\pi_{x}})(-\secx))\cdot\td(\T_{\pi_{x}})
=(\iota_{x}^{*}\ch(\T_{X})-\ch(\T_{\pi_{x}}))\cdot\ch(\mathcal{O}_{\ux}(-\secx))\cdot\td(\T_{\pi_{x}}).$$
La formule
$$\iota_{x}^{*}c_{1}(\X)=\ell(\secx+\pi_{x}^{*}c_{1}(\L_{x}))$$
ainsi que les formules
$$\iota_{x}^{*}c_{1}(\X)\cdot\secx=0$$
et
$$\secx^{2}=p^{*}c_{1}(\E_{x})\cdot\secx-c_{2}(\E_{x})=-\secx\cdot\pi_{x}^{*}c_{1}(\L_{x})  
\quad\text{(\voir\cite{Fu98} Remark 3.2.4)}$$ 
donnent
$$\secx\cdot \pi_{x}^{*}c_{1}(\L_{x})=\frac{1}{\ell^{2}}\iota_{x}^{*}c_{1}(\X)^{2}-\pi_{x}^{*}c_{1}(\L_{x})^{2}$$
et
$$\iota_{x}^{*}c_{1}(\X)\cdot \pi_{x}^{*}c_{1}(\L_{x})=\frac{1}{\ell^{2}}\iota_{x}^{*}c_{1}(\X)^{2}.$$
Le calcul donne (\voir\cite{Fu98} Example 3.2.3 et Example 3.2.4) :
\begin{eqnarray*}
\ch(\mathcal{O}_{\ux}(-\secx))_{\le 2} & = & 1-\secx+\frac{1}{2}{\secx}^{2} \\
 & = & 1-\secx+\frac{1}{2}(\pi_{x}^{*}c_{1}(\L_{x})^{2}-\frac{1}{\ell^{2}}\iota_{x}^{*}c_{1}(\X)^{2}),\\
\ch(\T_{\pi_{x}})_{\le 2} & = & 1+2\secx+\pi_{x}^{*}c_{1}(\L_{x})
+\frac{1}{2}(4{\secx}^{2}+\pi_{x}^{*}c_{1}(\L_{x})^{2}+4\secx\cdot\pi_{x}^{*}c_{1}(\L_{x}))\\
 & = & 1+2\secx+\pi_{x}^{*}c_{1}(\L_{x})+\frac{1}{2}\pi_{x}^{*}c_{1}(\L_{x})^{2},\\
\td(\T_{\pi_{x}})_{\le 2} & = & 1+\secx+\frac{1}{2}\pi_{x}^{*}c_{1}(\L_{x})
+\frac{1}{12}(4{\secx}^{2}+\pi_{x}^{*}c_{1}(\L_{x})^{2}+4\secx\cdot\pi_{x}^{*}c_{1}(\L_{x}))\\
 & = & 1+\secx+\frac{1}{2}\pi_{x}^{*}c_{1}(\L_{x})+\frac{1}{12}\pi_{x}^{*}c_{1}(\L_{x})^{2},
\end{eqnarray*}
et
$$\iota_{x}^{*}\ch(\T_{X})_{\le 2}=\iota_{x}^{*}(n+c_{1}(\X)+\frac{1}{2}(c_{1}(\X)^{2}-2c_{2}(\X))).$$
Finalement,
$$(\ch((\iota_{x}^{*}\T_{X}/\T_{\pi_{x}})(-\secx))\cdot\td(\T_{\pi_{x}}))_{2}=
\iota_{x}^{*}((\frac{\ell+1}{2\ell}-\frac{1}{\ell^{2} })c _{1}^{ 2}(\X)
-c_{2}(\X))+\frac{n-1}{12}\pi_{x}^{*}c_{1}(\L_{x})^{2}$$
et
\begin{eqnarray*}
c_{1}(\vx) & = & {\pi_{x}}_{*}(\ch((\iota_{x}^{*}\T_{X}/\T_{\pi_{x}})(-\secx))\cdot\td(\T_{\pi_{x}}))_{2}\\
 & = & {\pi_{x}}_{*}\iota_{x}^{*}((\frac{\ell+1}{2\ell}-\frac{1}{\ell^{2}})c_{1}^{2}(\X)-c_{2}(\X))
\in\textup{Pic}(\vx)\otimes\QQ.
\end{eqnarray*}
\end{proof}
\begin{rems}
Le calcul pr\'ec\'edent reste vrai si $x\in\X$ est
quelconque et si $\hx$ est une composante connexe de $\homxlibre\cap\hombir$.

L'extension 
$$0\longrightarrow\mathcal{O}_{\vx}\longrightarrow\E_{x}\longrightarrow\M_{x}\longrightarrow 0$$
est scind\'ee. Soit $\alpha$ sa classe dans $\textup{H}^{1}(\vx,\L_{x})$.
La restriction ${\pi_{x}}_{|\secx}$ de $\pi_{x}$ \`a $\secx$ est un isomorphisme de $\secx$ sur $\vx$. Nous noterons
encore $\L_{x}$ le fibr\'e ${\pi_{x}}_{|\secx}^{*}\L_{x}$ et $\alpha$ la classe
${\pi_{x}}_{|\secx}^{*}\alpha\in\textup{H}^{1}(\secx,\L_{x})$. 
\begin{lem}[\voir\cite{Dr04} Lemmes 3.2 et 3.3]
La classe du fibr\'e $\mathcal{O}_{\PP_{\vx}(\E_{x})}(1)\otimes \pi_{x}^{*}\L_{x}$ 
dans 
$$\textup{Pic}(2\secx)\simeq\textup{Pic}(\secx)\oplus\textup{H}^{1}(\secx,\L_{x})$$ 
est $(0,\alpha)$.
\end{lem}
\noindent Le fibr\'e $(\mathcal{O}_{\PP_{\vx}(\E_{x})}(-1)\otimes \pi_{x}^{*}\M_{x})^{\otimes\ell}$ est isomorphe au
fibr\'e $\iota_{x}^{*}\omega_{\X}$. La classe $\ell\alpha$ est donc nulle et $\alpha$ l'est aussi.
\end{rems}
\section{D\'emonstration des r\'esultats}
\begin{lem}\label{calculcycle}
Soit $\Y$ une vari\'et\'e alg\'ebrique connexe, projective et lisse sur $\CC$. Soit $\E$
un fibr\'e vectoriel de rang 2 sur $\Y$, extension d'un fibr\'e en droites
$\M$ par le fibr\'e trivial $\mathcal{O}_{\Y}$. Soient $\Z:=\PP_{\Y}(\E)$ et $p:\Z\longrightarrow\Y$ 
le morphisme naturel. Soit $\sigma\subset\Z$ la section de $p$ correspondant 
au quotient inversible $\M$ de $\E$.
Supposons qu'il existe une vari\'et\'e alg\'ebrique $\X$ projective et lisse sur $\CC$
et un morphisme $q:\Z\longrightarrow\X$ qui contracte $\sigma$ sur un point de $\X$.
Si $\beta$ est un cycle de codimension pure $k>0$ sur $\X$ alors
$$q_{*}(\Z)\cdot\beta=q_{*}(p^{*}\alpha\cdot(\sigma-p^{*} c_{1}(\M)))$$
o\`u $\alpha=p_{*}q^{*}\beta$ est un cycle de codimension $k-1$ sur $\Y$.
\end{lem}
\begin{proof}
L'anneau $\A(\Z)$ est un module libre sur $\A(\Y)$ engendr\'e
par la classe fondamentale de $\Z$ et la premi\`ere classe de Chern du fibr\'e tautologique 
(\voir\cite{Fu98} Theorem 3.3). 
Ici, la premi\`ere classe de
Chern du fibr\'e tautologique est $\sigma$ et
$$\sigma^{2}=p^{*}c_{1}(\E)\cdot\sigma-c_{2}(\E)=p^{*}c_{1}(\M)\cdot\sigma\quad\text{(\voir\cite{Fu98} Remark
3.2.4)}.$$ 
Le cycle $q^{*}\beta$ s'\'ecrit donc
$$q^{*}\beta=p^{*}\gamma\cdot\sigma+p^{*}\delta$$ 
avec $\gamma\in\A^{k-1}(\Y)$ et $\delta\in\A^{k}(\Y)$.
Le cycle $\gamma$ est donn\'e par la formule
$$p_{*}q^{*}\beta=\gamma\cdot p_{*}\sigma=\gamma.$$
Enfin, $\sigma$ \'etant contract\'ee par $q$ sur un point de $\X$ et le cyle $\beta$ \'etant de codimension
$>0$, on a 
$$0=q^{*}\beta\cdot\sigma=p^{*}(\gamma\cdot c_{1}(\M)+\delta)\cdot\sigma.$$
Les cycles $\gamma$ et $\delta$ v\'erifient donc la relation
$$\delta=-\gamma\cdot c_{1}(\M)$$
et
$$q_{*}(\Z)\cdot\beta=q_{*}q^{*}\beta=q_{*}(p^{*}\alpha\cdot(\sigma-p^{*} c_{1}(\M))).$$
\end{proof}
\begin{lem}\label{pointbase}
Si $\Y$ est une vari\'et\'e alg\'ebrique connexe, projective et lisse sur $\CC$ de dimension $n$ et si $\L$ est un
fibr\'e
en droites ample et engendr\'e par ses sections globales sur $\Y$ alors le
syst\`eme lin\'eaire $|\K_{\Y}+nc_{1}(\L)|$
est sans point base sauf si $(\Y,\L)\simeq(\PP^{n},\mathcal{O}_{\PP^{n}}(1))$. 
\end{lem}
\begin{proof}Le diviseur $\K_{\Y}+nc_{1}(\L)$ est num\'eriquement effectif si et seulement si
$(\Y,\L)\not\simeq(\PP^{n},\mathcal{O}_{\PP^{n}}(1))$ (\voir\cite{Fu87} Theorem 1). Supposons 
$(\Y,\L)\not\simeq(\PP^{n},\mathcal{O}_{\PP^{n}}(1))$ et montrons le r\'esultat
par r\'ecurence sur la dimension $n$ de $\Y$. 

Si $n=1$ alors le r\'esultat est imm\'ediat. Supposons $n\ge 2$. Soit $\H\in|\L|$ un diviseur g\'en\'eral. 
Le groupe $\textup{H}^{1}(\mathcal{O}_{\Y}(\K_{\Y}+(n-1)\H))$ est nul par le
th\'eor\`eme d'annulation de Kodaira.
L'application de restriction 
$$\textup{H}^{0}(\Y,\mathcal{O}_{\Y}(\K_{\Y}+n\H))\longrightarrow
\textup{H}^{0}(\H,\mathcal{O}_{\H}(\K_{\Y}+n\H))$$
est donc surjective. Le diviseur $(\K_{\Y}+n\H)_{|\H}=\K_{\H}+(n-1)\H_{|\H}$ est num\'eriquement effectif par la
formule d'adjonction. Si le syt\`eme lin\'eaire $|\K_{\H}+(n-1)\H_{|\H}|$ est sans
point base alors $|\K_{\Y}+n\H|$ l'est donc \'egalement.
\end{proof}
Les notations sont celles des paragraphes \ref{resultat}, \ref{tangent} et \ref{calcul}. Si $\hx'$ est une composante
irr\'eductible de $\hx$ nous noterons $\ux'$, $\vx'$, $\secx'$, $\iota_{x}'$, $\pi_{x}'$, $\M_{x}'$ et $\L_{x}'$ les
vari\'et\'es, morphismes et fibr\'es correspondants.
\begin{thm}\label{effectif}
Sous les hypoth\`eses \ref{setup}, si $x\in\X$ est g\'en\'eral 
et si $\ux'$ est une composante irr\'eductible de $\ux$
alors le cycle 
$${\iota_{x}'}_{*}(\ux')\cdot(c_{2}(\X)-\frac{\ell-1}{2\ell}c_{1}(\X)^{2})\in\A_{\ell-3}(\X)\otimes\QQ$$ 
est effectif.
\end{thm}
\begin{proof}
Si $x\in\X$ est g\'en\'eral alors les morphismes
param\'etr\'es par $\hx'$ sont libres (\voir Lemme \ref{libre}) et $\vx'$ est donc
lisse sur $\CC$ de dimesion $\ell-2$, o\`u $\ell:=c_{1}(\X)\cdot f_{*}\PP^{1}$ avec $[f]\in\hx'$ 
(\voir paragraphe \ref{tangent}). Si $x\in\X$ est toujours un point g\'en\'eral et si $[f]\in\hx'$ alors la courbe
$f(\PP^{1})$ est
immerg\'ee en $x$ (\voir\cite{Ke02} Theorem 3.3) : la restriction de la diff\'erentielle de l'application
$\iota_{x}'$ \`a $\T_{\pi_{x}'}\subset\T_{\ux'}$ 
d\'efinit un morphisme
$$\tau_{x}'\,:\,\vx'\longrightarrow\PP(\T_{\X,x}^{*})$$
qui \`a $[f]\in\vx'$ associe la droite $\textup{d}f_{0}(\T_{\PP^{1},0})$ de
$\T_{\X,x}$ (\voir\cite{Ke02} 3.2). Le fibr\'e en droites $\L_{x}'$ est isomorphe 
au fibr\'e ${\tau_{x}'}^{*}\mathcal{O}_{\PP(\T_{\X,x}^{*})}(1)$ et, puisque $\tau_{x}'$ est fini 
(\voir\cite{Ke02} Theorem 3.4), 
$\L_{x}'$ est ample et engendr\'e par ses sections globales. 

\medskip

Soit $\alpha:=-c_{1}(\vx')+(\ell-1)c_{1}(\L_{x}')\in\A_{\ell-3}(\vx')$.
Les formules (\voir paragraphe \ref{calcul})
$${\pi_{x}'}_{*}{\iota_{x}'}^{*}c_{1}(\X)^{2}=\ell^{2}c_{1}(\L_{x}')$$ 
et
$$c_{1}(\vx')={\pi_{x}'}_{*}{\iota_{x}'}^{*}((\frac{\ell+1}{2\ell}-\frac{1}{\ell^{2}})c_{1}(\X)^{2}-c_{2}(\X))$$
donnent
$${\pi_{x}'}_{*}{\iota_{x}'}^{*}(c_{2}(\X)-\frac{\ell-1}{2\ell}c_{1}(\X)^{2})=
\alpha$$
et, d'apr\`es le lemme \ref{calculcycle},
$${\iota_{x}'}_{*}(\ux')\cdot(c_{2}(\X)-\frac{\ell-1}{2\ell}c_{1}(\X)^{2})
={\iota_{x}'}_{*}({\pi_{x}'}^{*}\alpha\cdot(\secx'+{\pi_{x}'}^{*}c_{1}(\L_{x}'))).$$

\medskip

Le cycle $\alpha$ est effectif d'apr\`es le
lemme \ref{pointbase} et 
${\pi_{x}'}^{*}\alpha\cdot\secx'$ et $\alpha\cdot c_{1}(\L_{x}')$ le sont donc aussi. 
Le th\'eor\`eme est d\'emontr\'e.
\end{proof}
\begin{thm}\label{principal}
Sous les hypoth\`eses \ref{setup}, si $x\in\X$ est g\'en\'eral et si
$\ux'$ est une composante irr\'eductible de $\ux$
alors le cycle 
$${\iota_{x}'}_{*}(\ux')\cdot(c_{2}(\X)-(\frac{\ell-1}{2\ell}+\frac{1}{\ell^{2}})c_{1}(\X)^{2})\in\A_{\ell-3}(\X)
\otimes\QQ$$ 
est effectif sauf s'il existe un morphisme fini 
$\widehat{\X}\longrightarrow\X$ et une application rationnelle
$\varphi:\widehat{\X}\dashrightarrow\Z$ dont les fibres g\'en\'erales sont des
espaces projectifs sur $\CC$ de dimension $\ell-1$, tels que les courbes rationnelles consid\'er\'ees soient les
images dans $\X$ des droites contenues dans les fibres de $\varphi$.
\end{thm}
\begin{proof}
Soit maintenant $\alpha:=-c_{1}(\vx')+(\ell-2)c_{1}(\L_{x}')\in\A_{\ell-3}(\vx')$ le cycle tel que (\voir Lemme
\ref{calculcycle})
$${\iota_{x}'}_{*}(\ux')\cdot(c_{2}(\X)-(\frac{\ell-1}{2\ell}+\frac{1}{\ell^{2}})c_{1}(\X)^{2})
={\iota_{x}'}_{*}({\pi_{x}'}^{*}\alpha\cdot(\secx'+{\pi_{x}'}^{*}c_{1}(\L_{x}'))).$$
\noindent Le cycle $\alpha$ est effectif sauf si $(\vx',\L_{x}')\simeq(\PP^{\ell-2},\mathcal{O}_{\PP^{\ell-2}}(1))$
d'apr\`es le lemme \ref{pointbase}.
Si $\alpha$ est effectif alors ${\pi_{x}'}^{*}\alpha\cdot\secx'$ et $\alpha\cdot c_{1}(\L_{x}')$ le sont \'egalement et
le
cycle
$${\iota_{x}'}_{*}(\ux')\cdot(c_{2}(\X)-(\frac{\ell-1}{2\ell}+\frac{1}{\ell^{2}})c_{1}(\X)^{2})$$
est donc effectif dans ce cas.
Si $(\vx',\L_{x}')\simeq(\PP^{\ell-2},\mathcal{O}_{\PP^{\ell-2}}(1))$ alors
$${\iota_{x}'}_{*}(\ux')\cdot(c_{2}(\X)-\frac{\ell-1}{2\ell}c_{1}(\X)^{2})=0$$
et le cycle
$${\iota_{x}'}_{*}(\ux')\cdot(c_{2}(\X)-(\frac{\ell-1}{2\ell}+\frac{1}{\ell^{2}})c_{1}(\X)^{2})
=-{\iota_{x}'}_{*}({\pi_{x}'}^{*}c_{1}(\L_{x}')\cdot(\secx'+{\pi_{x}'}^{*}c_{1}(\L_{x}')))$$
n'est pas effectif sauf si $\ell=2$, auquel cas, il est nul.

\medskip

Supposons maintenant que l'ensemble $\Lambda$ des points de $\X$ tels que le cycle
$${\iota_{x}'}_{*}(\ux')\cdot(c_{2}(\X)-(\frac{\ell-1}{2\ell}+\frac{1}{\ell^{2}})c_{1}(\X)^{2})$$
ne soit pas effectif pour au moins une composante irr\'eductible $\ux'$ de $\ux$ soit dense dans $\X$. 
Si $x\in\Lambda$ alors, d'apr\`es la discussion qui pr\'ec\`ede, il existe une composante irr\'eductible $\vx'$ de $\vx$
telle que $(\vx',\L_{x}')\simeq(\PP^{\ell-2},\mathcal{O}_{\PP^{\ell-2}}(1))$. Supposons que
pour un point $x\in\X$ g\'en\'eral et toute composante connexe $\vx'$ de $\vx$, il existe des isomorphismes
$(\vx',\L_{x}')\simeq(\PP^{\ell-2},\mathcal{O}_{\PP^{\ell-2}}(1))$. 
Si $x\in\X$ est toujours g\'en\'eral et 
si $\vx'$ est une composante connexe  de $\vx$ alors $\tau_{x}'$ est une immersion ferm\'ee et $\tau_{x}'(\vx')$ est 
un sous-espace lin\'eaire de $\PP(\T_{\X,x}^{*})$. L'existence
d'un morphisme fini 
$\widehat{\X}\longrightarrow\X$ et d'une application rationnelle
$\varphi:\widehat{\X}\dashrightarrow\Z$ dont les fibres g\'en\'erales sont des
espaces projectifs sur $\CC$ de dimension $\ell-1$, tels que les courbes rationnelles consid\'er\'ees soient les
images dans $\X$ des droites contenues dans les fibres de $\varphi$ sont alors d\'emontr\'ees dans \cite{Ar04} (Theorem
3.1).

\medskip

Supposons donc $\Lambda$ dense dans $\X$ et montrons que pour un point $x\in\X$ g\'en\'eral et toute composante connexe
$\vx'$ de $\vx$, $(\vx',\L_{x}')\simeq(\PP^{\ell-2},\mathcal{O}_{\PP^{\ell-2}}(1))$.
Soit
$$\textup{c}:\hombirn\longrightarrow\textup{Chow}(\X)$$
l'application cycle (\voir\cite{Ko96} Corollary I 6.9), donn\'ee par
$$[f]\mapsto f_{*}\PP^{1},$$
o\`u $\hombirn$ est la normalisation de $\hombir$.
Soit $\overline{\V}$ la normalisation de l'adh\'erence de l'image de $\H$ dans $\textup{Chow}(\X)$ et soit
$\overline{\U}$ la normalisation du cycle universel sur $\overline{\V}$.
Le sous-ensemble
$\textup{c}(\H)$ est ouvert dans
son adh\'erence. 
Son image inverse dans $\overline{\V}$ est le quotient g\'eom\'etrique $\V$ de
$\H$ par $\groupe$ et son image inverse dans $\overline{\U}$
est le quotient g\'eom\'etrique $\U$ de $\PP^{1}\times\H$ par $\groupe$ (\voir\cite{Mo79} Lemma 9).
Soient $\bar{\pi}$ et $\bar{\iota}$

\centerline{
\xymatrix{
\overline{\U}\ar[d]^-{\bar{\pi}} \ar[r]^-{\bar{\iota}} & \X\\
\overline{\V}
}
}
\noindent les morphismes naturels. Le morphisme $\bar{\pi}$ est une fibration en droites projectives. 
Le point essentiel est que
l'image de $\overline{U}\setminus\U$ par $\bar{\iota}$ est un ferm\'e strict de $\X$
puisque le sch\'ema $\vx$ est propre sur $\CC$ pour $x\in\X$ g\'en\'eral.
Si $x\in\X$ est g\'en\'eral alors $\bar\iota^{-1}(x)\subset\U$ et
$\bar{\iota}$ est lisse au
dessus d'un ouvert non vide de $\X$ (\voir \cite{Ko96} II 3.5.3 et Lemme \ref{libre}).
La restriction de la diff\'erentielle de l'application
$\bar{\iota}$ \`a 
$\T_{\bar{\pi}}\subset\T_{\overline{\U}}$ d\'efinit une application rationnelle $\tau$

\centerline{
\xymatrix{
& \PP(\Omega_{\X}^{1})\ar[d]\\
 \overline{\U} \ar[r]^-{\bar{\iota}}\ar@{-->}@/^1.2pc/[ru]^-{\tau} & \X
}
}
\noindent bien d\'efinie au dessus d'un ouvert non vide de $\X$ (\voir \cite{Ke02} Theorem 3.3). 
Soit $\hat{\iota}:\overline{\U}\longrightarrow\widehat{\X}$ la factorisation de Stein de $\bar{\iota}$. Le morphisme
$\widehat{\X}\longrightarrow\X$ est \'etale au dessus d'un ouvert de $\X$.
Soit $\hat{\tau}$
l'application rationnelle d\'efinie par la restriction de la diff\'erentielle de l'application
$\hat{\iota}$ \`a $\T_{\bar{\pi}}\subset\T_{\overline{\U}}$

\centerline{
\xymatrix{
\PP(\Omega_{\X}^{1}) \ar[dd] & & \PP(\Omega_{\widehat{\X}}^{1}) \ar@{-->}[ll]\ar[dd] \\
 & \overline{\U} \ar[dr]^-{\hat{\iota}}\ar[dl]_-{\bar{\iota}}\ar@{-->}[ru]^-{\hat{\tau}}\ar@{-->}[ul]_-{\tau}\\
\X & & \widehat{\X} \ar[ll]
}
}
\noindent qui est bien d\'efinie au dessus d'un ouvert non vide de $\X$.
Il faut enfin remarquer que le morphisme
$\ux\supset\secx\simeq\vx\longrightarrow\bar{\iota}^{-1}(x)\subset\U$
est un isomorphisme si $x\in\X$ est g\'en\'eral.

L'ensemble $\Lambda$ est dense dans $\X$ par hypoth\`ese. L'ensemble des points $\hat x\in\widehat{\X}$
tels que
$$(\hat{\iota}^{-1}(\hat x),{\hat{\tau}^{*}\mathcal{O}_{\PP(\Omega_{\widehat{\X}}^{1})}(1)}_{|\hat{\iota}^{-1}(\hat x)})
\simeq(\PP^{\ell-2},\mathcal{O}_{\PP^{\ell-2}}(1))$$ 
est donc dense dans $\widehat{\X}$. Le faisceau
$\omega_{\overline{\U}/\widehat{\X}}\otimes\hat{\tau}^{*}\mathcal{O}_{\PP(\Omega_{\widehat{\X}}^{1})}(\ell-1)$
est plat au dessus d'un ouvert non vide de $\X$. Si $\hat x\in\widehat{\X}$ est un point
g\'en\'eral alors le faisceau 
 ${\omega_{\overline{\U}/\widehat{\X}}\otimes
 \hat{\tau}^{*}\mathcal{O}_{\PP(\Omega_{\widehat{\X}}^{1})}(\ell-1)}_{|\hat{\iota}^{-1}(\hat x)}$
est engendr\'e par ses sections globales (\voir Lemme
\ref{pointbase}) et trivial au dessus d'un point de $\Lambda$. L'ensemble des points de $\widehat{\X}$ o\`u l'espace
des sections de ce faisceau est de dimension 1 est en particulier ouvert dans
$\widehat{\X}$ (\voir \cite{Ha77} Theorem 12.8) et non vide par hypoth\`ese.
Le faisceau ${\omega_{\overline{\U}/\widehat{\X}}\otimes
 \hat{\tau}^{*}\mathcal{O}_{\PP(\Omega_{\widehat{\X}}^{1})}(\ell-1)}_{|\hat{\iota}^{-1}(\hat x)}$ est donc trivial pour
$\hat x\in\widehat{\X}$ g\'en\'eral et, d'apr\'es le crit\`ere de Kobayashi et Ochiai (\voir \cite{KO73}),
$$(\hat{\iota}^{-1}(\hat x),{\hat{\tau}^{*}\mathcal{O}_{\PP(\Omega_{\widehat{\X}}^{1})}(1)}_{|\hat{\iota}^{-1}(\hat x)})
\simeq(\PP^{\ell-2},\mathcal{O}_{\PP^{\ell-2}}(1)).$$
ce qui d\'emontre le r\'esultat annonc\'e plus haut.
\end{proof}
\begin{prop}\label{egalite}
Sous les hypoth\`eses \ref{setup}, le cycle
$${\iota_{x}}_{*}(\ux)\cdot(c_{2}(\X)-\frac{\ell-1}{2\ell}c_{1}(\X)^{2})\in\A_{\ell-3}(\X)\otimes\QQ$$ est nul
pour $x\in\X$ g\'en\'eral si et seulement s'il existe un morphisme fini 
$\widehat{\X}\longrightarrow\X$ et une application rationnelle
$\varphi:\widehat{\X}\dashrightarrow\Z$ dont les fibres g\'en\'erales sont des
espaces projectifs sur $\CC$ de dimension $\ell-1$, tels que les courbes rationnelles consid\'er\'ees soient les
images dans $\X$ des droites contenues dans les fibres de $\varphi$.
\end{prop}
\begin{proof}
Si le cycle
$${\iota_{x}}_{*}(\ux)\cdot(c_{2}(\X)-\frac{\ell-1}{2\ell}c_{1}(\X)^{2})\in\A_{\ell-3}(\X)\otimes\QQ$$ est nul
pour $x\in\X$ g\'en\'eral alors  
$$(\vx',\L_{x}')\simeq(\PP^{\ell-2},\mathcal{O}_{\PP^{\ell-2}}(1))$$ 
pour toute composante connexe $\vx'$ de $\vx$ (\voir Th\'eor\`eme \ref{effectif}) et
le Theorem 3.1 de \cite{Ar04} donne la conclusion cherch\'ee.

Inversement, supposons qu'il existe un morphisme fini $\nu:\widehat{\X}\longrightarrow\X$ et 
une application rationnelle $\varphi:\widehat{\X}\dashrightarrow\Z$ dont les fibres 
g\'en\'erales sont des espaces projectifs sur $\CC$ de dimension $\ell-1$, tels que les courbes
rationnelles consid\'er\'ees soient les images dans $\X$ des droites contenues dans les fibres de $\varphi$. Soit
$\F\simeq\PP^{\ell-1}$ une fibre g\'en\'erale de $\varphi$ et $j:\F\hookrightarrow\widehat{\X}$ l'immersion ferm\'ee de
$\F$ dans $\widehat{\X}$. La formule de projection
donne 
\begin{eqnarray*}
\F\cdot\nu^{*}(c_{2}(\X)-\frac{\ell-1}{2\ell}c_{1}(\X)^{2}) & = &
j_{*}(c_{2}(\F)-\frac{\ell-1}{2\ell}c_{1}(\F)^{2})\\
 & = & 0
\end{eqnarray*}
et
\begin{eqnarray*}
{\iota_{x}}_{*}(\ux)\cdot(c_{2}(\X)-\frac{\ell-1}{2\ell}c_{1}(\X)^{2}) & = & 
\nu_{*}(\F)\cdot(c_{2}(\X)-\frac{\ell-1}{2\ell}c_{1}(\X)^{2})\\
& = & \nu_{*}(\F\cdot\nu^{*}(c_{2}(\X)-\frac{\ell-1}{2\ell}c_{1}(\X)^{2}))\\
& = & 0.
\end{eqnarray*}
\end{proof}
\begin{rem}
Supposons que la diff\'erentielle de $\tau_{x}'$ s'annule en codimension 1 et
soit $\dx'\subset\vx'$ le diviseur de $\vx'$ tel que l'application induite
$$\T_{\vx'}\longrightarrow{\tau_{x}'}^{*}\T_{\PP(\T_{\X,x}^{*})}(-\dx')$$
soit non nulle en codimension 1.
Son conoyau $\Q$ est localement libre en codimension 1. Le fibr\'e
${\tau_{x}'}^{*}\T_{\PP(\T_{\X,x}^{*})}\otimes\M_{x}'$
\'etant engendr\'e par ses sections globales, le faisceau $\Q\otimes\mathcal{O}_{\vx'}(\dx')\otimes\M_{x}'$ l'est
\'egalement.
\begin{lem}\label{cycleeffectif}
Soit $\Y$ une $\CC$-vari\'et\'e r\'eguli\`ere de dimension $n\ge 1$. Si $\G$ est un faisceau coh\'erent,
localement libre en codimension 1 et engendr\'e par
ses sections globales alors le cycle $c_{1}(\G)\in\A(\Y)$ est effectif.
\end{lem}
\begin{proof}
Le faisceau $\G$ est localement libre sur un ouvert $\Y_{0}\subset\Y$ dont le compl\'ementaire est de codimension au
moins 2 dans $\Y$ et $\G_{|\Y_{0}}$ est \'egalement engendr\'e par ses sections globales.
Notons $i$ l'inclusion de $\Y_{0}$ dans $\Y$. L'application $i^{*}:\A_{n-1}(\Y)\longrightarrow\A_{n-1}(\Y_{0})$
\'etant un isomorphisme, il suffit de traiter le cas o\`u $\G$ est localement libre et engendr\'e par ses sections
globales, qui est imm\'ediat.
\end{proof}
\noindent Soit $\alpha:=c_{1}(\Q\otimes\mathcal{O}_{\vx'}(\dx')\otimes\M_{x}')\in\A_{\ell-3}(\vx')$. Le cycle $\alpha$
est
effectif d'apr\`es le lemme \ref{cycleeffectif},
\begin{eqnarray*}
c_{1}(\vx') & = & c_{1}(\L_{x}')+(n-1)(c_{1}(\L_{x}')-\dx')-c_{1}(\Q)\\
& = & c_{1}(\L_{x}')+(n-1)(c_{1}(\L_{x}')-\dx')+(n-\ell+1)(-c_{1}(\L_{x}')+\dx')-\alpha\\
& = & (\ell-1)c_{1}(\L_{x}')-(\ell-2)\dx'-\alpha
\end{eqnarray*}
et
$${\iota_{x}'}_{*}(\ux')\cdot(c_{2}(\X)-\frac{\ell-1}{2\ell}c_{1}(\X)^{2})=
\frac{1}{2\ell}{\iota_{x}'}_{*}({\pi_{x}'}^{*}((\ell-2)\dx'+\alpha)\cdot(\secx'+{\pi_{x}'}^{*}c_{1}(\L_{x}'))).$$
La description g\'eom\'etrique du support de $\dx'$ est donn\'ee par le r\'esultat suivant.
\begin{lem}[\voir\cite{Ar04} Proposition 2.7]Le  morphisme $\tau_{x}'$ est une immersion en $[f]\in\vx'$ si et seulement
si $\displaystyle{f^{*}\T_{\X}\simeq\mathcal{O}_{\PP^{1}}(2)\oplus\mathcal{O}_{\PP^{1}}(1)^{\oplus \ell-2}\oplus
\mathcal{O}_{\PP^{1}}^{\oplus n-\ell+1}}.$
\end{lem} 
\end{rem}


\begin{thebibliography}{abcdefg}
\bibitem[Ar04]{Ar04}C. Araujo, \textit{Rational curves of minimal degree and characterizations of $\PP^{n}$}, preprint 
math.AG/04010584.
\bibitem[CO75]{CO75}B. Chen, K. Ogiue, \textit{Some characterizations of complex space forms in terms of Chern
classes}, Quart. J. Math. Oxford \textbf{26} (1975), 459-464.
\bibitem[De01]{De01}O. Debarre, \textit{Higher-dimensional algebraic geometry}, Universitext, Springer-Verlag, 2001.
\bibitem[DPS94]{DPS94}J.-P. Demailly, T. Peternell, M. Schneider, \textit{Compact complex manifolds with numerically
effective tangent bundles}, J. of Algebraic Geometry \textbf{3} (1994), 295-345.
\bibitem[Dr04]{Dr04}S. Druel, \textit{Caractérisation de l'espace projectif}, Manuscripta Math. 
\textbf{115} (2004), 19-30.
\bibitem[Fu87]{Fu87}T. Fujita, \textit{On polarized manifolds whose adjoint bundles are not semipositive}, Adv. Stud.
Pure Math. \textbf{10} (1987), 167-178.
\bibitem[Fu98]{Fu98}W. Fulton, \textit{Intersection theory}, Second edition, Ergebnisse der Mathematik und ihre
Grenzgebiete \textbf{3}, Springer-Verlag, 1998.
\bibitem[Ha77]{Ha77}R. Hartshorne, \textit{Algebraic Geometry}, Graduate Texts in 
Mathematics \textbf{52}, Springer-Verlag, 1977.
\bibitem[Ke02]{Ke02}S. Kebekus, \textit{Families of singular rational curves}, J. of Algebraic Geometry \textbf{11}
(2002), 245-256.
\bibitem[KO73]{KO73}S. Kobayashi, T. Ochiai, \textit{Characterization of
complex projective spaces and hyperquadrics}, J. Math. Kyoto Univ. \textbf{13} (1973),
31-47.
\bibitem[Ko96]{Ko96}J. Koll\'ar, \textit{Rational curves on algebraic varieties},
Ergebnisse der Mathematik und ihre Grenzgebiete \textbf{32}, Springer-Verlag, 1996.
\bibitem[L\"u82]{Lu82}M. L\"ubke, \textit{Chernklassen von Hermite-Einstein Vectorb\"undeln}, Math. Ann. \textbf{260}
(1982), 133-141.
\bibitem[Mo79]{Mo79}S. Mori, \textit{Projective manifolds with ample
tangent bundles},  Ann. of Math. \textbf {110} (1979), 593--606.
\bibitem[Ti02]{Ti02}G. Tian, \textit{Canonical metrics in K\"ahler geometry}, Lect. in Math. ETH, Birkh\"auser (2000).
\end{thebibliography}
\end{document}